\documentclass[a4paper]{amsart}
\usepackage{a4,pstricks,bbm}
\usepackage{amsmath,amssymb,latexsym}
\usepackage{epsf,psfrag,epsfig}
\usepackage[english]{babel}
\usepackage[latin1]{inputenc}

 \address{Instituto de Ciencias Matem\'{a}ticas (CSIC-UAM-UC3M-UCM), 28049 Madrid, Spain}

  \email{juanjo.rue@icmat.es}

\newtheorem{thm}{Theorem}

\newcommand{\vu}{\mathbf{m}}
\newcommand{\findem}{\hfill $\square$\\ \medskip}

\newcommand{\cA}{\mathcal{A}}

\title{On polynomial representation functions for multivariate linear forms}
\author{Juanjo Ru\'e}
\thanks{The author is supported by a JAE-DOC grant from the Junta para la Ampliaci\'on de Estudios (CSIC), jointly
financed by the FSE, by the MTM2011-22851 grant, (Spain) and the ICMAT Severo Ochoa project SEV-2011-0087 (Spain).}

\begin{document}

\maketitle

\begin{flushright}
\emph{Yahya Ould Hamidoune, in Memoriam.}
\end{flushright}

\begin{abstract}
Given an infinite sequence of positive integers $\cA$, we prove that for every non-negative integer $k$ the number of solutions of the equation $n=a_1+\dots+a_k$, $a_1,\,\dots, a_k\in \cA$, is not constant for $n$ large enough. This result is a corollary of our main theorem, which partially answers a question of S\'ark\"ozy and S\'os on representation functions for multivariate linear forms. Additionally, we obtain an Erd\H{o}s-Fuchs type result for a wide variety of representation functions.
\end{abstract}

\section{Introduction}
Let $\cA$ be an infinite sequence of positive integers. Denote by $r(n,\cA)$ the number of solutions of the equation $n=a_1+a_2$, where $a_1,a_2\in \cA$. In~\cite{ErdTuran} the authors found, by means of analytic arguments, that $r(n,\cA)$ cannot be constant for $n$ large enough. As it is shown in~\cite{Dirac}, and elementary argument also exists: it is obvious that $r(n,\cA)$ is odd when $n = 2a,\, a\in\cA $, and even otherwise. So it is not possible that $r(n,\cA)$ is constant for $n$ large enough. This idea can be easily generalized when we consider the number of solutions of the equation $n=a_1+\dots+a_p$, where $a_1,\dots, a_p\in \cA$ and $p$ is a prime number: if $a\in \cA$, the number of representations of $pa$ is congruent to 1 modulo $p$, while the number of representations of $pa+1$ is congruent to $0$ modulo $p$. As $a$ can be chosen as big as desired, a contradiction is obtained if we suppose that the representation function is constant for $n$ large enough. However, the argument fails when we consider a composite modulo, and it does not seem that the argument could be extended in the general case using elementary tools.

These problems are particular cases of a question posted by S\'ark\"ozy and S\'os~\cite{SarSos}: given a multivariate linear form $k_1x_1+\dots+k_r x_r$, consider the number of solutions of the equation $n=k_1 a_1+\dots+k_r a_r$, where $a_1,\dots,a_r\in \cA$. For which multivariate linear forms the number of solutions could be constant for $n$ large enough? For bivariate linear forms with positive coefficients the problem is completely solved: when we deal with the bivariate linear form $x_1+k x_2$, $k>1$, Moser~\cite{Moser} showed that there exists a set $\cA$ such that the number of solutions of the equation $n=a_1+ka_2$ where $a_1,\,a_2\in \cA$ is constant and equal to 1 (see also~\cite{Shev} for additional properties of these sequences of numbers). Recently, Cilleruelo and Ru\'e~\cite{CillRue} proved that for $k_1$ and $k_2$ satisfying $1<k_1<k_2$ and $\mathrm{gcd}(k_1,k_2)=1$ the number of solutions of the equation $n=k_1 a_1+k_2 a_2$ where $a_1,\,a_2\in \cA$ is not constant for $n$ large enough.

In this paper we found an answer to the question posed by S\'ark\"ozy and S\'os for multivariate linear forms in several cases: let $0<k_1<k_2<\dots<k_r$ be a finite sequence of positive integers and consider the \emph{configuration} $\vu=\{(k_1,m_1),\dots, (k_r,m_r)\}$, for $m_1,\dots, m_r>0$. Each value $m_i$ is the \emph{multiplicity} of $k_i$, and the \emph{degree} of $\vu$ is $\mathrm{gcd}(m_1,\dots, m_r)$. Given a configuration $\vu$, we consider the associated multivariate linear form $k_1 \left(x_{1,1}+\dots+x_{1,m_1}\right)+\dots+k_r\left( x_{r,1}+\dots+ x_{r,m_r}\right)$. Given a sequence of positive integers $\cA$, the \emph{representation function} of $n$ with respect to $\vu$ is the number of different solutions of the equation
\begin{equation}\label{eq:repr-func}
n=k_1 \left(a_{1,1}+\dots+a_{1,m_1}\right)+\dots+k_r\left( a_{r,1}+\dots+ a_{r,m_r}\right)
\end{equation}
where $a_{i,j}\in \cA$. We denote this value by $r_{\vu}(n,\cA)$. Our first theorem deals with representation functions which are polynomials:
\begin{thm}\label{thm:main}
Let $\cA$ be an infinite sequence of positive integers, and $\vu$ a configuration of degree $s$. Then, no polynomial of degree smaller than $s-1$ can coincide with the function $r_{\vu}(n,\cA)$ for $n$ large enough.
\end{thm}
As a trivial consequence, Theorem~\ref{thm:main} solves the problem posted by S\'ark\"ozy and S\'os for multivariate linear forms associated to configurations whose degree is greater than 1. Observe also that the exponent $s-1$ in Theorem~\ref{thm:main} cannot be improved in general. To see this, let $\cA=\mathbb{N}$ and consider the multivariate linear form $x_{1,1}+\dots+x_{1,s}$. Then, the associated configuration is $\textbf{m}=\{(1,s)\}$ and the number of representations of $n$ is $r_{\textbf{m}}(n,\cA)=\binom{n+1}{s-1}$, which is a polynomial of degree $s-1$.

In this paper we also deal with a question related to Erd\H{o}s-Fuchs Theorem, motivated by the lattice point problem: if we write $S(n)=\left|\{(x,y)\in \mathbb{Z}^2\, : x^2+y^2 \leq n \}\right|-\pi n$, an estimate by Hardy and Landau of $S(n)$ states that  $S(n)=o\left(n^{1/4}\left(\log n\right)^{1/4}\right)$ cannot hold. Erd\H{o}s-Fuchs Theorem states a similar result for arbitrary sets: let $\cA$ be an infinite sequence of positive integers and $\varepsilon>0$. Then
\begin{equation}\label{eq:erdos-fuchs-1}
\sum_{j=1}^{n}\left(r(j,\cA)-c\right)=O\left(n^{1/4-\varepsilon}\right)
\end{equation}
cannot hold for any positive constant $c$ (recall that $r(j,\cA)$ is the number of solutions of the equation $a_1+a_2=j$ with $a_1,\,a_2\in \cA$). The precise bound obtained by Erd\H{o}s and Fuchs is $o\left(n^{1/4}\left(\log n\right)^{-1/2}\right)$ instead of $O\left(n^{1/4-\varepsilon}\right)$. Further improvements of this theorem~\cite{montgomery-vaughan, newman-simplified} show that a similar result also holds when the right hand side of~\eqref{eq:erdos-fuchs-1} is equal to $o\left(n^{1/4}\right)$, and in fact this bound is the best possible~\cite{ruzsa}. Some results in this direction has been deduced for several summands~\cite{min-tang,horvath-manysummands}, and also when considering the sum of different sequences~\cite{horvath-improvement,sarkozy}.

In this paper we get an Erd\H{o}s-Fuchs type result for representation functions associated to configurations with degree greater than $1$. In particular, we have the following theorem:
\begin{thm}\label{thm:erdos-fuchs}
Let $\cA$ be an infinite sequence of positive integers, let $\vu$ be a configuration of degree $s>1$ and $\varepsilon>0$. Then,
\begin{equation*}\label{eq:erdos-fuchs}
\sum_{j=1}^{n}\left(r_{\vu}(j,\cA)-c\right)=O\left(n^{1/4-\varepsilon}\right)
\end{equation*}
cannot hold for any positive constant $c$.
\end{thm}
In particular, Theorem $2$ implies Theorem $1$ in the case of polynomials of degree $0$. The method used to get the proof of Theorem~\ref{thm:erdos-fuchs} is divided in two cases, depending on the parity of $s$. Our arguments generalize the one that appear in~\cite{newman-simplified}. Observe that our bound is rougher than the $o\left(n^{1/4}\right)$ that is known for other representation functions. However, our argument is shorter and simpler compared with the one needed to get this exponent.
\\
\paragraph{\textbf{Plan of the paper:}} in Section~\ref{sec:tools} we introduce the necessary background in order to deal with the problem, namely the use of generating functions in order to codify the problem, Carlson's Theorem for power series with radius of convergence $1$. We prove Theorem~\ref{thm:main} in Section~\ref{sec:main}, and Theorem~\ref{thm:erdos-fuchs} in Section~\ref{sec:erdos-fuchs}.


\section{Tools} \label{sec:tools}

\paragraph{\textbf{Generating functions:}} we codify all the enumerative information of the problem using generating functions: given a set $\cA$ of non-negative integers we define the formal power series
$$f_{\cA}(z)=\sum_{a\in \cA} z^a.$$
This series is called the \emph{generating function} associated to the sequence $\cA$.  Since its coefficients are either $0$ or $1$, this formal power series is either a polynomial (if $\cA$ is finite) or has a singularity at $z=1$ (if $\cA$ is infinite). In the second case, the Taylor expansion of $f_{\cA}(z)$ around $z=0$ has radius of convergence equals to $1$, hence all its singularities have modulo greater or equal to $1$.

The combinatorial problem can be translated in the language of generating functions in the following way. Let $\cA$ be a sequence of non-negative integers and let $\vu=\{(k_1,m_1),\dots, (k_r,m_r)\}$. Then,
\begin{eqnarray*} \label{eq:main-gf}
\left(f_\cA\left(z^{k_1}\right)\right)^{m_1}\dots \left(f_\cA\left(z^{k_r}\right)\right)^{m_r}&=&\sum_{a_{i,j}\in \cA}z^{k_1 \left(a_{1,1}+\dots+a_{1,m_1}\right)+\dots+k_r\left(a_{r,1}+\dots+ a_{r,m_r}\right)}\nonumber\\
&=&\sum_{n=0}^{\infty}r_{\vu}(n,\cA)z^{n}.
\end{eqnarray*}

\paragraph{\textbf{Carlson's Theorem:}} further details on this result could be found in~\cite{remmert}. We denote by $\mathbb{E}$ the disk $\{u\in \mathbb{C}: |u|<1\}$. Carlson's Theorem assert the following dichotomy:

\begin{thm}[Carlson's Theorem] \label{thm:carlson}
Let $f(z) = \sum_{n=0}^{\infty} a_nz^n$  be a power series with integer coefficients and radius of convergence $R = 1$. Then either $\mathbb{E}$  is the domain of holomorphy of $f(z)$ or $f(z)$ can be extended to a rational function of the form $S(z)/(1 - z^m)^n$, where $S(z)\in \mathbb{Z}[z]$  and $m, n \in  \mathbb{N}$.
\end{thm}

\section{Proof of Theorem~\ref{thm:main}}\label{sec:main}

Without loss of generality, we may assume that $0\in \cA$. We suppose that such a sequence $\cA$ exists, and we argue by contradiction. The case $s=1$ is trivial, so we may assume that $s>1$. We assume that $r_{\vu}(n,\cA)$ is equal to a polynomial of degree at most $s-2$, namely $q(n)=\sum_{i=0}^{d}q_i n^i$, with $d<s-1$, for $n$ large enough, say $n>N$. Using the generating function terminology,
\begin{equation}\label{eq:polyn}
\sum_{n=0}^{\infty}r_{\vu}(n,\cA)z^{n}=T(z)+\sum_{n>N}^{\infty}q(n)z^n=T(z)+\sum_{i=0}^{d}q_i\sum_{n>N}^{\infty}n^{i}z^n=T_0(z)+\sum_{i=0}^{d}q_i\sum_{n=0}^{\infty}n^{i}z^n,
\end{equation}
where $T(z),\,T_0(z)$ are polynomials with degree  $\leq N$. Each term of the form $\sum_{n=0}^{\infty}n^{i}z^n$ can be written as $\frac{Q_i(z)}{(1-z)^{i+1}}$, where $Q_i(z)$ is a polynomial in $z$ such that $Q_i(1)\neq 0$. Hence, we can write Expression~\eqref{eq:polyn} in the form
$$T_0(z)+\sum_{i=0}^{d}q_i\frac{Q_i(z)}{(1-z)^{i+1}}=\frac{Q(z)}{(1-z)^{d+1}},$$
where $Q(z)$ is a polynomial which satisfies that $Q(1)\neq 0$. Using now Equation~\eqref{eq:repr-func}, the generating function $f_\cA(z)$ satisfies the relation
\begin{equation}\label{eq:eq-main}
\left(f_\cA\left(z^{k_1}\right)\right)^{m_1}\dots \left(f_\cA\left(z^{k_r}\right)\right)^{m_r}=\frac{Q(z)}{(1-z)^{d+1}}.
\end{equation}
Observe that  $Q(0)=T(0)=r_{\vu}(0,\cA)=1$. As the degree of  $\vu$ is equal to $s$, Equation~\eqref{eq:eq-main} can be written in the form
\begin{equation*}\label{eq:eq-main2}
\left(\left(f_\cA\left(z^{k_1}\right)\right)^{m_1/s}\dots \left(f_\cA\left(z^{k_r}\right)\right)^{m_r/s}\right)^s=\frac{Q(z)}{(1-z)^{d+1}},
\end{equation*}
where each quotient $m_i/s$ is a non-negative integer (recall that each $m_i>0$).  Note that the Taylor development of $\left(f_\cA\left(z^{k_1}\right)\right)^{m_1/s}\dots \left(f_\cA\left(z^{k_r}\right)\right)^{m_r/s}$ has non-negative integer coefficients and radius of convergence equals to $1$ (as each term has a singularity at $z=1$). Carlson's Theorem asserts then that two situations may happen: $\left(f_\cA\left(z^{k_1}\right)\right)^{m_1/s}\dots \left(f_\cA\left(z^{k_r}\right)\right)^{m_r/s}$ is either a rational function or it has $\mathbb{E}$ as a domain of holomorphy. It is obvious that the second situation may not happen. Let us assume the first condition, namely
$$\left(f_\cA\left(z^{k_1}\right)\right)^{m_1/s}\dots \left(f_\cA\left(z^{k_r}\right)\right)^{m_r/s}= \frac{S(z)}{\left(1-z^m\right)^{n}},$$
where $S(z)\in \mathbb{Z}[z]$, and consequently
$$\left(\frac{S(z)}{\left(1-z^m\right)^{n}}\right)^s=\frac{Q(z)}{(1-z)^{d+1}}.$$
This relation could be written as the equality $S(z)^s(1-z)^{d+1}=Q(z)\left(1-z^m\right)^{ns}$. As $Q(1)\neq 0$, observe then that all the roots of $Q(z)$ have degree multiple of $s$, namely $Q(z)=p(z)^s$ for a certain polynomial $p(z)$. Hence $S(z)^s(1-z)^{d+1}=p(z)^s\left(1-z^m\right)^{ns}$ and $s$ divides $d+1<(s-1)+1=s$, which is a contradiction. \findem

\section{Proof of Theorem~\ref{thm:erdos-fuchs}}\label{sec:erdos-fuchs}

We write $\sum_{j=1}^{n}\left(r_{\vu}(n,\cA)-c\right)=a_n$. We assume that $a_n=O\left(n^{1/4-\varepsilon}\right)$. Using the generating function methodology, this condition can be written as
\begin{equation*}\label{eq:erdos-fuchs-gf}
\frac{1}{1-z}\left(f_\cA\left(z^{k_1}\right)\right)^{m_1}\dots \left(f_\cA\left(z^{k_r}\right)\right)^{m_r}=\frac{c}{\left(1-z\right)^2}+\sum_{n=0}^{\infty}a_n z^n.
\end{equation*}
As $\gcd(m_1,\dots,m_r)=s>1$, we write $F_\cA(z)=\left(f_\cA\left(z^{k_1}\right)\right)^{m_1/s}\dots \left(f_\cA\left(z^{k_r}\right)\right)^{m_r/s}$ in the previous equation, getting
\begin{equation}\label{eq:erdos-fuchs-main}
F_\cA(z)^s=\frac{c}{1-z}+(1-z)\sum_{n=0}^{\infty}a_n z^n.
\end{equation}
We write $m_i'=m_i/s$. Observe that the $n$th Taylor coefficient of $F_\cA(z)$ is the number of solutions of the equation
$$n=k_1\left(a_{1,1}+\dots+ a_{1,m_1'}\right)+\dots+k_r\left(a_{r,1}+\dots+ a_{r,m_r'}\right),\, a_{i,j}\in \cA.$$
Consequently, considering the configuration $\vu'=\{(k_1,m_1'),\dots, (k_r,m_r')\}$, we have $F_{\cA}(z)=\sum_{n=0}^{\infty}r_{\vu'}(n,\cA)z^n$.

Define $h_{M}(z)=1+z+\dots+z^{M-1}=\frac{1-z^M}{1-z}$. We start multiplying Equation~\eqref{eq:erdos-fuchs-main} by $h_M(z)^2$ , getting the equality
$$F_\cA(z)^s h_M(z)^2=\frac{c}{1-z}h_M(z)^2+h_M(z)^2(1-z)\sum_{n=0}^{\infty}a_n z^n,$$
which gives the inequality
\begin{equation}\label{eq:s=2-1}
\left|F_\cA(z)\right|^s\left|h_M(z)\right|^2 \leq \frac{cM^2}{|1-z|}+2\left|h_M(z)\sum_{n=0}^{\infty}a_n z^n\right|.
\end{equation}
Let $r$ be a positive real number slightly smaller than $1$ (later we will define it properly). The strategy of the proof is based on integrating~\eqref{eq:s=2-1} along the circle $\mathbb{S}_r=\{z\in \mathbb{C}:\, |z|=r\}$ in order to get bounds for both the left and right hand side. We distinguish two situations depending on whether $s$ is even or not.

\subsection{The even case}\label{subsect:2k}
Assume first that $s$ is even, namely $s=2k$. Starting with the left hand side, we write $F_\cA(z)^k h_M(z)=\sum_{n=0}^{\infty}b_n z^n$. Observe that coefficients $b_n$ are positive integers, so $b_n^2\geq b_n$ for each $n$. We proceed applying Parseval's Theorem; the following integral is considered with respect to the normalized arc length $\frac{|dz|}{2\pi r}$:
$$\frac{1}{2\pi r}\int_{\mathbb{S}_r} \left|F_\cA(z)^{k} h_M(z)\right|^2 |dz|=\sum_{n=0}^{\infty} b_n^2 r^{2n} \geq \sum_{n=0}^{\infty} b_n r^{2n}   = F_\cA\left(r^2\right)^k h_M\left(r^2\right) \geq F_\cA\left(r^2\right)^k M r^{2M} .$$

Let us get a bound for $F_\cA\left(r^2\right)$. As we are assuming that $a_n=O\left(n^{1/4-\varepsilon}\right)$, as $z\rightarrow 1^-$, $F_\cA(z)^{2k} \sim \frac{c}{1-z}$. Then, for $r$ close enough to $z=1$, there exists a constant $C_1$ such that
$$F_{\cA}\left(r^2\right)^{2k}=\frac{c}{1-r^2}+(1-r^2)\sum_{n=0}^{\infty}a_n r^{2n} \geq \frac{C_1}{1-r^2}.$$
We conclude with the following bound:
\begin{equation}\label{eq:bound1}
\frac{1}{2\pi r}\int_{\mathbb{S}_r} \left|F_\cA(z)^k h_M(z)\right|^2 |dz|\geq C_2\frac{M}{\left(1-r^2\right)^{1/2}}r^{2M},
\end{equation}
where $C_2=C_1^{1/2}$.

\subsection{The odd case}\label{subsect:2k+1}
In order to get a similar bound in the case where $s$ is an odd integer, namely $s=2k+1$, we apply H\"older's inequality over the left hand side of Equation~\eqref{eq:bound1}. Recall that for $u,v$ complex-valued functions, this inequality states that
\begin{equation}\label{eq:holder}
\|uv\|_1\leq \|u\|_p \|v\|_q,
\end{equation}
where $1\leq p,q< \infty$ and $\frac{1}{p}+\frac{1}{q}=1$. Choosing $p=\frac{2k+1}{2k}$, $q=2k+1$, $u=F_{\cA}(z)^{2k}h_M(z)^{\frac{4k}{2k+1}}$ and $v=h_M(z)^{\frac{2}{2k+1}}$ in~\eqref{eq:holder} and simplifying conveniently the exponents we get the inequality
%
%
\begin{equation}\label{eq:holder1}
\frac{1}{2\pi r}\int_{\mathbb{S}_r} \left|F_{\cA}(z)\right|^{2k+1}\left|h_M(z)\right|^2 |dz| \geq \frac{\left(\frac{1}{2\pi r}\int_{\mathbb{S}_r}|F_{\cA}(z)|^{2k} \left|h_M(z)\right|^2 |dz|\right)^{\frac{2k+1}{2k}}}{\left(\frac{1}{2\pi r}\int_{\mathbb{S}_r}\left|h_M(z)\right|^2 |dz|\right)^{\frac{1}{2k}}}
\end{equation}

We get bounds for both the numerator and the denominator of the right hand side of~\eqref{eq:holder1} by applying Parseval's Theorem. More concretely, for the numerator we apply an argument similar to the one used to get bound~\eqref{eq:bound1}, getting
\begin{eqnarray*}
\left(\frac{1}{2\pi r}\int_{\mathbb{S}_r}|F_{\cA}(z)|^{2k} \left|h_M(z)\right|^2 |dz|\right)^{\frac{2k+1}{2k}}&\geq& M^{\frac{2k+1}{2k}} r^{2M \frac{2k+1}{2k}}\left(F_\cA\left(r^2\right)^k\right)^{\frac{2k+1}{2k}}\\
&=&M^{\frac{2k+1}{2k}} r^{2M \frac{2k+1}{2k}} F_\cA\left(r^2\right)^{\frac{2k+1}{2}}\geq C_3\frac{M^{\frac{2k+1}{2k}} r^{2M \frac{2k+1}{2k}}}{(1-r^2)^{1/2}}  ,
\end{eqnarray*}
for a certain constant $C_3$, and for the denominator we integrate directly (using Parseval's Theorem)
$$\left(\frac{1}{2\pi r}\int_{\mathbb{S}_r}\left|h_M(z)\right|^2 |dz|\right)^{\frac{1}{2k}}= h_M\left(r^2\right)^{\frac{1}{2k}} \leq M^{\frac{1}{2k}}. $$
Consequently, we deduce that
\begin{equation}\label{eq:bound2}
\frac{1}{2\pi r}\int_{\mathbb{S}_r} \left|F_{\cA}(z)\right|^{2k+1}\left|h_M(z)\right|^2 |dz| \geq \frac{C_3\frac{M^{\frac{2k+1}{2k}} r^{2M \frac{2k+1}{2k}}}{(1-r^2)^{1/2}}}{M^{\frac{1}{2k}}}=C_3 \frac{M}{(1-r^2)^{1/2}}r^{2M \frac{2k+1}{2k}}.
\end{equation}
\subsection{The final argument} We continue with the right hand side of Equation~\eqref{eq:s=2-1}. Applying Parseval's Theorem on the first summand
\begin{eqnarray}\label{eq:log-1}
\frac{1}{2\pi r}\int_{\mathbb{S}_r}\frac{cM^2}{|1-z|} |dz|&=& c M^2\sum_{n=0}^{\infty}\frac{1}{4^{2n}} \binom{2n}{n}^2 r^{2n} \leq c M^2\left(1+\sum_{n=1}^{\infty}\frac{1}{n} r^{2n}\right) \\ &=& cM^2\left(1+\log\left(\frac{1}{1-r^2}\right)\right)< 2cM^2\log\left(\frac{1}{1-r^2}\right) ,\nonumber
\end{eqnarray}
and for the second summand, invoking Cauchy-Schwarz inequality we obtain
\begin{eqnarray}\label{eq:log-2}
\frac{1}{2\pi r}\int_{\mathbb{S}_r} \left|h_M(z)\sum_{n=0}^{\infty}a_n z^n\right| |dz|&\leq & \frac{1}{2\pi r} \left(\int_{\mathbb{S}_r} \left|h_M(z)\right|^2 |dz| \right)^{1/2} \left(\int_{\mathbb{S}_r} \left|\sum_{n=0}^{\infty}a_n z^n\right|^2 |dz|\right)^{1/2}\nonumber\\
&=& \left(\sum_{n=0}^{M-1}r^{2n}\right)^{1/2} \left(\sum_{n=0}^{\infty} a_n^2 r^{2n}\right)^{1/2} \leq M^{1/2} \left(\sum_{n=0}^{\infty} a_n^2 r^{2n}\right)^{1/2}.
\end{eqnarray}
In order to get a bound for the sum in~\eqref{eq:log-2}, we use that $a_n=O(n^{1/4-\varepsilon})$, hence
$$\sum_{n=0}^{\infty} a_n^2 r^{2n} = O\left(\sum_{n=0}^{\infty} n^{1/2-2\varepsilon} r^{2n}\right)=O\left(\frac{1}{(1-r^2)^{3/2-2\varepsilon}}\right).$$
Resuming, joining bounds \eqref{eq:log-1}, \eqref{eq:log-2} we have obtained
\begin{equation}\label{eq:dreta}
\frac{1}{2\pi r}\int_{\mathbb{S}_r}\frac{cM^2}{|1-z|} |dz|+\frac{1}{2\pi r}\int_{\mathbb{S}_r} \left|h_M(z)\sum_{n=0}^{\infty}a_n z^n\right| |dz| < 2cM^2\log\left(\frac{1}{1-r^2}\right)+C'\frac{M^{1/2}}{(1-r^2)^{3/4-\varepsilon}},
\end{equation}
where $C'$ is a constant.

For $M$ large enough, write $r^2=1-M^{-(2+8\varepsilon)}$. Observe that $1>r^{2M}=\left(1-M^{-1}\right)^{M} > \frac{1}{4}$ for all $M$, hence bounds in~\eqref{eq:bound1} and~\eqref{eq:bound2} are asymptotically equal (up to a constant term). Note also that that $\left(1-r^2\right)^{-1}=M^{2+8\varepsilon}$. Using this change of variables, bounds in Equation~\eqref{eq:bound1} and~\eqref{eq:bound2} are the same, so it is not necessary to distinguish between them. We write $C$ to denote the constant $C_2$ or $C_3$ (that is, either $C=C_2$ or $C=C_3$ depending on whether $s$ is even or odd, respectively). Hence, substituting these asymptotic bounds in Equation~\eqref{eq:dreta} we get
$$C M^{2+4\varepsilon}<4c(1+4\varepsilon) M^2\log\left(M\right)+C'M^{2+4\varepsilon-8\varepsilon^2},$$
or equivalently, for a certain constant $C$
\begin{equation}\label{eq:final}
C <4c(1+4\varepsilon) M^{-4\varepsilon}\log\left(M\right)+C'M^{-8\varepsilon^2}.
\end{equation}
Making $M\rightarrow \infty$ the right hand side of \eqref{eq:final} tends to $0$, and this contradiction finishes the proof of Theorem 2. \findem

\paragraph{\textbf{Acknowlegments:}} Javier Cilleruelo is greatly thanked for inspiring discussions and useful comments in order to improve the presentation of this paper. Maksym Radziwill is kindly appreciated for pointing Carlson Theorem in order to simplify the proof of Theorem~\ref{thm:main}. The author also thanks Oriol Serra and Boris Bukh for further suggestions related to the proof of Theorem~\ref{thm:erdos-fuchs}. Finally, the author also thank the referee, whose suggestions and advices helped to improve the presentation of the article.

{\small \bibliography{Juanjo-bib}}

\begin{thebibliography}{10}

\bibitem{CillRue}
{\sc Cilleruelo, J., and Ru\'e, J.}
\newblock On a question of {S}\'ark\"ozy and {S}\'os on bilinear forms.
\newblock {\em Bulletin of the London Mathematical Society 41}, 2 (2009),
  274--280.

\bibitem{Dirac}
{\sc Dirac, G.}
\newblock Note on a problem in additive number theory.
\newblock {\em Journal of the London Mathematical Society 26\/} (1951),
  312--313.

\bibitem{ErdTuran}
{\sc Erd\H{o}s, P., and Tur\'an, P.}
\newblock On a problem of {S}idon in additive number theory, and on some
  related problems.
\newblock {\em Journal of the London Mathematical Society 16\/} (1941),
  212--215.

\bibitem{horvath-manysummands}
{\sc Horv{\'a}th, G.}
\newblock On a theorem of {E}rd{\H o}s and {F}uchs.
\newblock {\em Acta Arith. 103}, 4 (2002), 321--328.

\bibitem{horvath-improvement}
{\sc Horv{\'a}th, G.}
\newblock An improvement of an extension of a theorem of {E}rd{\H o}s and
  {F}uchs.
\newblock {\em Acta Math. Hungar. 104}, 1-2 (2004), 27--37.

\bibitem{montgomery-vaughan}
{\sc Montgomery, H.~L., and Vaughan, R.~C.}
\newblock On the {E}rd{\H o}s-{F}uchs theorems.
\newblock In {\em A tribute to {P}aul {E}rd{\H o}s}. Cambridge Univ. Press,
  Cambridge, 1990, pp.~331--338.

\bibitem{Moser}
{\sc Moser, L.}
\newblock An application of generating series.
\newblock {\em Mathematics Magazine 1}, 35 (1962), 37--38.

\bibitem{newman-simplified}
{\sc Newman, D.~J.}
\newblock A simplified proof of the {E}rd{\H o}s-{F}uchs theorem.
\newblock {\em Proc. Amer. Math. Soc. 75}, 2 (1979), 209--210.

\bibitem{remmert}
{\sc Remmert, R.}
\newblock {\em Classical topics in complex function theory}, vol.~172 of {\em
  Graduate Texts in Mathematics}.
\newblock Springer, New York, 1998.

\bibitem{ruzsa}
{\sc Ruzsa, I.~Z.}
\newblock A converse to a theorem of {E}rd{\H o}s and {F}uchs.
\newblock {\em J. Number Theory 62}, 2 (1997), 397--402.

\bibitem{sarkozy}
{\sc S{\'a}rk{\"o}zy, A.}
\newblock On a theorem of {E}rd{\H o}s and {F}uchs.
\newblock {\em Acta Arith. 37\/} (1980), 333--338.

\bibitem{SarSos}
{\sc S{\'a}rk{\"o}zy, A., and S{\'o}s, V.~T.}
\newblock On additive representation functions.
\newblock In {\em The mathematics of {P}aul {E}rd{\H o}s, {I}}, vol.~13 of {\em
  Algorithms Combin.} Springer, Berlin, 1997, pp.~129--150.

\bibitem{Shev}
{\sc Shevelev, V.}
\newblock On unique additive representations of positive integers and some
  close problems.
\newblock {\em Available on-line at arXiv:0811.0290\/}.

\bibitem{min-tang}
{\sc Tang, M.}
\newblock On a generalization of a theorem of {E}rd{\H o}s and {F}uchs.
\newblock {\em Discrete Math. 309}, 21 (2009), 6288--6293.

\end{thebibliography}
\bibliographystyle{acm}

\end{document}